\documentclass[a4paper]{article}
\usepackage{amsmath}
\usepackage{amsfonts}
\usepackage{amsthm}
\usepackage{amssymb}
\usepackage[toc,page]{appendix}
\usepackage{array}
\usepackage{enumitem}
\usepackage{float}
\usepackage{setspace}
\usepackage{graphicx}
\usepackage{leftidx}
\usepackage{arydshln}
\usepackage{tikz-cd} 
\usepackage{hyperref}
\usepackage{mathrsfs}
\usepackage{etoolbox}
\usepackage{makeidx}
\usepackage{tikz}
\usepackage{soul}
\usepackage{mathtools}
\usetikzlibrary{decorations.pathmorphing}

\patchcmd{\thebibliography}{\section*{\refname}}{}{}{}

\newtheorem{Theorem}{Theorem}[section]

\newtheorem{Definition}[Theorem]{Definition}
\newtheorem{Meta-Definition}[Theorem]{Meta-Definition}

\newtheorem{Remark}[Theorem]{Remark}

\newtheorem{Conjecture}[Theorem]{Conjecture}
\newtheorem{Statement}[Theorem]{Statement}

\DeclareMathOperator{\re}{Re}

\DeclareMathOperator{\id}{id}

\DeclareMathOperator{\Aff}{\mathcal{A}ff}

\DeclareMathOperator{\STAR}{STAR}

\DeclareMathOperator{\Small}{small}
\DeclareMathOperator{\Larget}{large}

\title{On the monograph "Finiteness Theorems for limit cycles" and a special case of alternant cycles}

\author{Melvin Yeung}

\date{}

\begin{document}
	
	\maketitle
	
	\begin{abstract}
		We provide evidence that the approach of \cite{ilyashenkoFiniteness} to the proof of Dulac's theorem has a gap. Although the asymptotics of \cite{ilyashenkoFiniteness} capture far more than the asymptotics of Dulac, we prove that the arguments for why the asymptotics in \cite{ilyashenkoFiniteness} are not themselves oscillatory is insufficient. We give an explicit counterexample and we draw confines to which Ilyashenko's result may be restricted in order to keep the validity.
	\end{abstract}
	
	\section{Introduction}
	
	Dulac's theorem, asserted in \cite{DulacCycleslimites} with a heavily flawed proof, states that the number of limit cycles of any given polynomial vector field on the plane is finite. It is a soft version of the second part of the 16th problem of Hilbert, asking for a (sharp) upper bound on the number of limit cycles in terms of the polynomial degree.\\
	
	Let us sketch a very brief history, firmly restricting to Dulac's conjecture and not widening the history to the 16th problem of Hilbert. For a more general overview of the 16th problem of Hilbert, including a wide variety of related problems, see \cite{Ilyashenko2002CentennialHO}.\\
	
	In Dulac's original paper, \cite{DulacCycleslimites}, several reductions were carried out, the primary reduction being based on Poincar\'e-Bendixson's theorem. After compactifying the phase plane to a sphere, he noted that an unbounded sequence of limit cycles should accumulate to either an equilibrium, a periodic orbit, or a graphic (so a finite union of equilibria with orbits connecting them).\\
	
	Because of compactness one only has to prove that limit cycles cannot accumulate onto any of the sets above.  Ruling out the possibility of accumulating onto a periodic orbit is not difficult at all. On the other hand both equilibria and graphics contain an incredibly rich structure which can not immediately be studied. In order to study them anyway we use a version of a blowup at a point desingularizing them to elementary graphics which we will call polycycles. The price to pay here is that our polycycles are defined on 2-dimensional real analytic manifolds instead of just neatly on a sphere. It is worth noting that the fact that this process can be done in a finite number of steps was only proven much later, see e.g. \cite{DUMORTIER197753}.\\ 
	
	Polycycles in this context are formed by a finite number of equilibria, all of them being hyperbolic or semi-hyperbolic (so with one nonzero eigenvalue) of the vector field, all of them being connected with regular heteroclinic or homoclinic connections (without additional singularities on the connections themselves). Without loss of generality we will only look at a single side of a polycycle at a time, so for example a cuspidal loop will be studied on the inside and the outside separately.\\
	
	Dulac then proceeded to study compositions of transition maps defined near hyperbolic saddles or semi-hyperbolic saddles.  Next Dulac made an additional, crucial, reduction to so-called balanced cycles.  By grouping the semi-hyperbolic passages in groups of two, he showed the existence of an asymptotic expansion for these Dulac maps in terms of $(x,\log x)$ and proved a group property of such maps.  The mistake he made was to take the quasi-analyticity for granted: he inferred the triviality of the map from the triviality of its asymptotic expansion.\\
	
	Ilyashenko in \cite{IlyashenkoMemoir} produced a clever counter example, clearly showing why Dulac's arguments failed and additionally he showed that Dulac's theorem is valid for hyperbolic polycycles, i.e. polycycles with only hyperbolic equilibria. It is a corner stone in this story and up to date the only result that has not been questioned.\\
	
	The author would like to thank Dmitry Novikov for helping him understand the quasianalyticity arguments of \cite{ilyashenkoFiniteness} and for the feedback on this article.\\
	
	The author is very grateful to Daniel Panazzolo for organizing a workshop for the author to present his findings and for the feedback on this article.\\
	
	The author thanks his promotor Peter de Maesschalck for introducing him to this problem and for the extensive feedback on this article.\\
	
	The author is thankful to Yulij S. Ilyashenko for providing feedback on this article. 
	
	\subsection{The known case and the main difficulty}
	
	\subsubsection{Hyperbolic polycycles}

	With the notion of almost regular map, Ilyashenko has successfully treated the hyperbolic case, see e.g. the introduction of \cite{ilyashenkoFiniteness}: almost regular maps are maps that are analytically continuable to a so-called standard quadratic domain, and that are asymptotic to a Dulac-asymptotic series inside this domain.\\
	
	The size of the domain of the analytic continuation, together with the presence of the asymptotic series gives the quasianalyticity property: an almost regular map that is asymptotic to the identity map up to exponential accuracy, \emph{is} the identity map.\\  
	
	Furthermore, almost regular maps form a group with respect to the composition so any first return map defined near some hyperbolic polycycle is actually an almost regular map.  So either it is asymptotic to the identity (and then it is the identity, so near the polycycle all orbits are periodic), or it is not, and then by looking at the most dominant term in the expansion it is easy to conclude that there are no periodic orbits accumulating to the polycycle.
	
	\subsubsection{Adding more depth}
	
	Semi-hyperbolic singularities present in the polycycle are much more delicate due to the exponential flatness of the Dulac map of a semi-hyperbolic saddle. The standard example being the following:
	
	\[E(x) = e^{-\frac{1}{x}}\]
	
	which is the Dulac map of:
	
	\[\begin{cases*}
		\dot{x} = x^{2}\\
		\dot{y} = -y
	\end{cases*}\]
	
	from a section $\left\{y = \frac{1}{e}\right\}$ transverse to the hyperbolic separatrix to a section $\{x = 1\}$ transverse to the center separatrix, phase plane restricted to the positive quadrant. The point is that conjugation of a function $f$ with $E$ allows us to make flatter and flatter asymptotics, which may no longer be picked up if your asymptotics is not sensitive enough. Let us write a function:
	
	\[f(x) = x + x^{2}f_{1}(x).\]
	
	Then:
	
	\begin{align*}
		E^{-1}(f(E(x)) &= \frac{-1}{\ln\left(E(x) + E(x)^{2}f_{1}(E(x))\right)}\\ 
		&= -\frac{1}{-\frac{1}{x} + \ln(1 + E(x)f_{1}(E(x)))}\\ 
		&= x\frac{1}{1 - x\ln(1 + E(x)f_{1}(E(x)))}.
	\end{align*}
	
	This function is identity $+$ flat terms, thus having Dulac series $x$ but nonetheless not being equal to $x$. More complicated examples are certainly possible, think for example of a cycle giving first return map:
	
	\[E^{-1}\circ f \circ E \circ g \circ E^{-1} \circ E^{-1} \circ h \circ E \circ E.\]
	
	Note that by glueing as in \cite[\S0.3 C]{ilyashenkoFiniteness} such a polycycle should definitely exist on some 2-dimensional real analytic manifold.\\
	
	Now we see an interplay of three levels of asymptotics: level 0 coming from $g$, level 1 coming from $E^{-1}\circ f \circ E$ and level 2 coming from $E^{-2} \circ h \circ E^{2}$.  It should however by clear by now why so-called unbalanced cycles are easily dealt with: either they are flat or inverse to flat. While in specific situations a clever choice of transversal or looking at the difference map with respect to two transversals instead of the return map can resolve some issues, in general no level of flatness can be ignored.	
	
	\subsection{Current state of the proofs}
	
	\subsubsection{Commonalities}
	
	Before diving into the differences of the two proofs in \cite{ecalleconstructiveproof} and the main subject of the article today \cite{ilyashenkoFiniteness} we will first talk in broad terms about what both approaches share.\\
	
	The consensus holds that the spirit of the proof of Dulac was correct, what one should do is to find \textit{sufficiently accurate} asymptotics for first return maps of polycycles in order to determine the whole map.\\
	
	Note that it is in the first place not clear what asymptotics should mean in this context, for example the previous section makes clear that it should be necessary to know what it means for a map to be asymptotic to:
	
	\[\left(\sum_{n \geq 0}n!x^{n}\right) + e^{-\frac{1}{x}}.\]
	
	The consensus between \cite{ecalleconstructiveproof} and \cite{ilyashenkoFiniteness} is also that (in the real analytic context) \textit{formal} information should be enough to determine these asymptotics in the following sense: Suppose given a semihyperbolic saddle, then there are formal coordinates (i.e. a coordinate transform in terms of diverging power series) such that it is orbitally equivalent to a semihyperbolic saddle of the form:
	
	\[\begin{cases*}
		\dot{x} = \frac{x^{k + 1}}{1 + ax^{k}}\\
		\dot{y} = -y
	\end{cases*}\]
	
	they agree that what one should do is compose the asymptotics coming from the formal normalization maps, the transit maps of the normal forms above, the usual Dulac series for the hyperbolic saddles and the Taylor series for the transit maps.\\
	
	This gives an expression called a transseries, a formal object that is able to represent asymptotic information ``at different depths". The notion of transseries is defined explicitly in \cite{ecalleconstructiveproof} and contained implicitly in \cite{ilyashenkoFiniteness} under the name STAR-series. A STAR-series is not exactly the same as a transseries, it is a mix of both the asymptotic information above and the way to construct a sensible asymptotic series for the first return map of a polycycle.\\
	
	Transseries has been the center of quite vibrant research including a recent result linking this to O-minimality, see \cite{ADHTransseries}, offering a road forward to the 16th problem of Hilbert once the roadblock of understanding these proofs have passed. It is also worth noting that transseries have been used in contexts more explicitly related to the problem of Dulac, for example in \cite{GalalKaiserSpeisegger2020}. An older result in the same vein \cite{KaiserRolinSpeissegger2009} in fact gives uniform bounds for limit cycles close to a polycycle with only nonresonant hyperbolic saddles.\\
	
	The strategy in \cite{ecalleconstructiveproof} and \cite{ilyashenkoFiniteness} is to recast the transserial formal asymptotics in a form which says something meaningful about the functions they are describing, i.e. return maps of polycycles. There should be at minimum two properties of these asymptotics:
	
	\begin{enumerate}
		\item{If they are nonzero they should give a leading term.}
		\item{If they are zero they should only describe the zero function, i.e. the quasianalyticity property.}
	\end{enumerate}
	
	\subsubsection{Ecalle's approach}
	
	In their essence, the asymptotic objects of Ecalle are simple, they are the transseries we spoke of. Thus simply assuming compatibility with the addition, we can just take out the leading term in any nonzero transseries by standard theory, see e.g. \cite{ADHTransseries}, giving leading terms for nonzero transseries.\\
	
	The part that remains difficult is the claims of quasianalyticity. From what the author understands of the matter the claim is that \cite{ecalleconstructiveproof} provides an injection from a certain class of transseries, called accelero-summable transseries, into germs of functions near infinity (the polycycle return maps being brought there under the coordinate transform $\frac{1}{x}$). It is worth noting that transseries in \cite{ecalleconstructiveproof} are not exactly the same as in \cite{ADHTransseries} but are instead more general, see \cite[p 149, (4.1.59)]{ecalleconstructiveproof} for a transseries not included in \cite{ADHTransseries}.\\
	
	Again, from what the author understands, in \cite[Chapter 5]{ecalleconstructiveproof}, something called cryptolinear formulas are introduced which split up transseries into blocks which are individually amenable to a generalized version of a Borel transform, which after a change of coordinates is amenable to a Laplace transform which needs to be `accelerated back into the right time' and then everything can be summed back up, this process should then be accelero-summation, giving the injection above.\\
	
	While at the time of writing, Ecalle's proof, \cite{ecalleconstructiveproof} is still not being debated, the proof has not been digested by the community either, and although many of Ecalle's claims have meanwhile been properly proved by others using his methods, see e.g. \cite{Costin_2019}, \cite{SauzinSemihyperbolic}, his proof of Dulac's theorem is still waiting for an accessible and detailed version.\\
	
	Considering the amount of effort put into this problem, it is hard to believe that in all this time no one has been able to draw a detailed proof from the texts provided by Ecalle. In some sense, this calls for a following conjecture:
	
	\begin{Conjecture}
		\cite{ecalleconstructiveproof} can be worked out to a fully detailed proof without additional difficulties.
	\end{Conjecture}
	
	In particular looking at \cite[page 157-158]{ecalleconstructiveproof} the case of convergent normalizing maps is mentioned and in this case the claim is that this can be `naively summed' as described in \cite[pp. 142-143]{ecalleconstructiveproof} which seems a much more conventional definition. Perhaps it would be interesting to work out this case first.
	
	\subsubsection{Ilyashenko's approach}
	
	Conversely the proof in \cite{ilyashenkoFiniteness} has a very good answer for why the quasianalyticity should hold, in an extension of his proof in the hyperbolic case the answer is simply that the return maps of polycycle can be decomposed into function which have small Stokes phenomena on a large domain, the domains being large enough and the Stokes phenomena being small enough that it implies quasi-analyticity.\\
	
	Seeing as this is the main subject of the current article we will refrain from technical detail in this part.\\
	
	Recently, Ilyashenko has come to realize a flaw in his proof himself, and he communicated openly on this matter during a conference at the Fields institute in 2021. The flaw communicated by Ilyashenko lies in the fact that while any semihyperbolic saddle can be analytically conjugated to a system of the form:
	
	\[\begin{cases*}
		\dot{x} = \frac{x^{k + 1}}{1 + ax^{k}}\\
		\dot{y} = -y
	\end{cases*}\]
	
	this conjugation does not necessarily map the real axis to the real axis. In particular in the cases where the domain is large enough to complete his quasianalyticity arguments, it is known that the conjugation maps in general do not map the real axis to the real axis. In fact, as is shown in \cite{MartinetRamisSemihyperbolic}, reality is the exception rather than the norm.\\
	
	This results in the fact that the asymptotics in \cite{ilyashenkoFiniteness} are complex valued. In order to try to keep closer to the ordered field of the reals and get leading terms anyway \cite{ilyashenkoFiniteness} essentially uses the construction that for any complex analytic function $f$ on the real axis the function:
	
	\[z \mapsto \frac{f(z) + \overline{f(\bar{z})}}{2}\]
	
	is a real analytic function which is $\re(f)$ on the real axis. The problem Ilyashenko has communicated is that this can not be done while preserving the needed quasi-analyticity properties.\\ 
	
	In 2020, the author took on the challenge to attempt to both salvage the problem and to make Ilyashenko's proof more accessible. (In fact, the challenge first consisted solely of the second part, but the communication of the flaw in the proof added an extra level.)  He succeeded in giving an easier presentation, immediately restricting to the case where the real axis gets normalized to the real axis.\\
	
	Nevertheless the author was able to verify that the quasi-analyticity argument of \cite{ilyashenkoFiniteness} holds even in full generality for all polycycles, it was the existence of a leading term which was an issue.\\
	
	Shortly after the author communicated his findings at a workshop in Mulhouse, France, he found an additional mistake of Ilyashenko's proof (which corresponds to a gap in his own version of the proof). This mistake pertains to the existence of a leading term coming from the asymptotics given in \cite{ilyashenkoFiniteness}. In some sense unlike the conjecture before which was possibly born of insufficient understanding, we have been able to understand enough of \cite{ilyashenkoFiniteness} to prove in this article that:
	
	\begin{Statement}\label{StatementBase}
		The proof for Dulac's theorem, provided by Ilyashenko in \cite{ilyashenkoFiniteness}, cannot be worked out to full
		detail without additional difficulties.\\
		
		In particular here the difficulties lie in the existence of a leading term for a nonzero asymptotic series.
	\end{Statement}
	
	It is a new element that casts a shadow on all proofs.  What is worse is that the presence of two independent proofs of Dulac's result blocked any further development on the matter, except for some spurious cases.  Partial results, for example for polycycles with exactly two semi-hyperbolic points (see \cite{IlyashenkoSeparatrixLunes}), have therefore not had the amount of attention that they should have had.\\
	
	Even so, just like \cite{ilyashenkoFiniteness} and \cite{ecalleconstructiveproof} saw it fit to continue the approach of \cite{DulacCycleslimites} and try to fill in the gaps because of the inherent value of the arguments of \cite{DulacCycleslimites}, the author feels that \cite{ilyashenkoFiniteness} is valuable. Specifically the entire approach laid out in \cite{ilyashenkoFiniteness} working towards the Additive Decomposition Theorem \cite[p 73]{ilyashenkoFiniteness} using Shift Lemma's should be worked out to completion, filling in as much as possible the gap laid out in this article. This roughly corresponds to the arguments in Sections \ref{SectGeom} and \ref{SectFormal} of this article.\\
	
	In a personal communication with Ilyashenko, he announced a work in progress, \cite{IlyashenkoRewritten}, aiming to cover the so-called depth 1 or alternant case (more on this later) with a new method that avoids the arguments of \cite[\S 4.10]{ilyashenkoFiniteness} as a possible answer to the reality issues mentioned earlier. Instead choosing an approach based on normal forms, Ilyashenko is aiming at addressing the extra difficulties announced in Statement \ref{StatementBase} and bypassing Statement \ref{StatementCounterexample} as well, in the alternant case.
	
	\section{Goals}
	
	\begin{Statement}\label{StatementCounterexample}
		At the core of the existence of a leading term for a Dulac map we have identified and isolated the following incorrect claim at the end of \cite[p 198]{ilyashenkoFiniteness} that (calculation error corrected in the last term):
		
		\[a_{1}'a_{2} - a_{1}a_{2}' + \mathbf{e}'a_{1}a_{2} \in \mathscr{K}^{m, r}_{\mathbb{R}}\]
		
		which is at the core of the existence of a leading term for a Dulac map in the way the proof goes in \cite{ilyashenkoFiniteness}.
	\end{Statement}
	
	The statement above may seem banal at a glance, yet here we make the claim that it is at the core of the existence of a leading term, so we will first dedicate most of the text to explaining why this is at the core of the proof.\\
			
	We will do this restriction to a much simpler class of polycycles than the general case and simplifying the arguments of \cite{ilyashenkoFiniteness} to this context.\\
	
	\begin{Statement}
		The problem isolated in Statement \ref{StatementCounterexample} already occurs in polycycles satisfying the following conditions:
		
		\begin{enumerate}
			\item{All equilibria in the polycycle are for some $k \geq 1$, locally real analytically orbitally equivalent to the following system at the origin (in forward or backward time):
				
				\[\begin{cases*}
					\dot{x} = x^{k + 1}\\
					\dot{y} = -y.
				\end{cases*}\]}
			\item{The equilibria `alternate' in the following sense, suppose that one equilibrium is traversed going from the hyperbolic side to the center side, i.e. in the system above coming from $y > 0, x$ small positive going to to $x > 0, y$ small positive, then the next orbit in forward time will be traversed going from the center side to the hyperbolic side.}
		\end{enumerate}
		
		This is a special case of so-called alternant polycycles, term coined by Ilyashenko, where all semihyperbolic saddles/saddle nodes in the polycycle satisfy the second condition.\\
		
		Therefore we will call these polycycles simple alternant polycycles.
	\end{Statement}
	
	\begin{Remark}
		Alternant polycycles are exactly the depth 1 case that the article in preparation \cite{IlyashenkoRewritten} will cover.
	\end{Remark}
		
	In the simplification of these arguments we will deviate from \cite{ilyashenkoFiniteness} in the sense that we will talk in terms of minimal conditions for the arguments of \cite{ilyashenkoFiniteness} to work, i.e. for the central result, the additive decomposition theorem (\cite[p 73]{ilyashenkoFiniteness}) to work rather than the explicit classes that are given in \cite{ilyashenkoFiniteness}.\\
		
	We want to show that complications as in point $1$ are absolutely impossible to avoid in the methods of \cite{ilyashenkoFiniteness}.\\
	
	Then we will present our counterexample in section \ref{SectCounterExample}.\\
		
	Next we will restrict a bit further, nonetheless to a case not covered by the well-known arguments of \cite{ilyashenkoFiniteness}, i.e. the hyperbolic case, and we will show that these arguments do work in this case.
	
	\begin{Statement}\label{StatementPos}
		For simple alternating polycycles where every equilibrium is real analytically orbitally equivalent to:
		
		\[\begin{cases*}
			\dot{x} = x^{2}\\
			\dot{y} = -y.
		\end{cases*}\]
		
		The difficulty with statement \ref{StatementCounterexample} disappears and the proof of \cite{ilyashenkoFiniteness} is valid.
	\end{Statement}
	
	Let us describe the structure of the article. Sketching out the proof of \cite{ilyashenkoFiniteness} will take in the next three sections coinciding roughly with the three large parts of the proof of \cite{ilyashenkoFiniteness}:
	
	\begin{enumerate}
		\item{In Section \ref{SectGeom} we will prove the Structural Theorem which can be roughly considered all the `dynamics' of \cite{ilyashenkoFiniteness}, after this section we might as well forget we are talking about polycycles.}
		\item{In Section \ref{SectFormal} we will present all of the formal manipulations in \cite{ilyashenkoFiniteness}, roughly corresponding to \cite[\S 1.11]{ilyashenkoFiniteness}.}
		\item{In Section \ref{SectProblem} we will define what in our case would be the $\mathscr{K}^{m, r}_{\mathbb{R}}$ of Statement \ref{StatementCounterexample} and we will explain why Statement \ref{StatementCounterexample} is crucial in the existence of a leading term, essentially expositing the arguments in \cite[\S 4.10]{ilyashenkoFiniteness}.}
	\end{enumerate}
	
	Then we will provide the counterexample to Statement \ref{StatementCounterexample} in Section \ref{SectCounterExample}, followed by proving Statement \ref{StatementPos} in Section \ref{SectPosres}, even though most of the work happens in Section \ref{SectCounterExample}.
	
	\section{Splitting up the return map and the logarithmic chart}\label{SectGeom}
	
	The first step in the proof following \cite{ilyashenkoFiniteness} is to split up the return map into transit maps. By definition of simple alternant polycycle near each equilibrium $E_{i}$ there exists real analytic sections $\Sigma_{i1}, \Sigma_{i2}$ transversally crossing the separatrices for which the transition map $T_{i}$ is either given by:
	
	\[z \mapsto e^{-\frac{1}{z^{k_{i}}}}\]
	
	for some $k_{i} \geq 1$ or its inverse:
	
	\[z \mapsto \left(\frac{1}{-\ln\left(z\right)}\right)^{\frac{1}{k_{i}}}.\]
	
	Suppose the equilibria ordered from $1$ to $N$ in order of which we encounter them. We will assume the transversal $\Sigma$ relative to which we take the Dulac map to be $\Sigma_{N2}$.
	
	\begin{figure}[H]
		\centering
		\includegraphics[scale = 1]{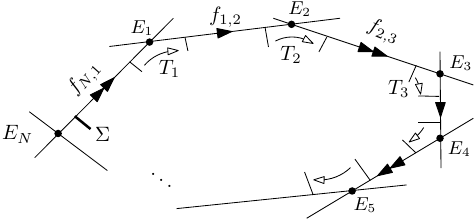}
		\caption{Decomposition of a first return map of a simple alternant polycycle.  The arrows on the separatrices are double for a hyperbolic separatrix and single for a center separatrix.}
	\end{figure}
	
	Let $f_{i, j}$ be the analytic flow box map from $\Sigma_{i2}$ to $\Sigma_{j1}$, then the return map $\Delta$ is given by
	
	\[\Delta = T_{N} \circ f_{N - 1, N}\circ T_{N - 1} \circ \cdots \circ T_{1} \circ f_{N, 1}.\]
	
	Now we may assume without loss of generality (by rechoosing $\Sigma$) we start on a hyperbolic separatrix, i.e. $T_{1} = e^{-\frac{1}{z^{k_{1}}}}$.\\
	
	The next step is to introduce the so-called logarithmic chart, basically because the maps $e^{-\frac{1}{z}}$ and $\frac{1}{-\ln(z)}$ are very complicated, the second one even implying a branch cut of the logarithm. In order to simplify these we introduce the logarithmic chart:
	
	\[\zeta = -\ln(z).\]
	
	For a function $g$ in the usual $z$ chart we denote its counterpart in the logarithmic chart by $g^{\log}$, more concretely:
	
	\[g^{\log}(\zeta) = -\ln(g(e^{-\zeta})).\]
	
	Now note that $\left(z \mapsto e^{-\frac{1}{z}}\right)^{\log} = \zeta \mapsto e^{\zeta}$ and $\left(z \mapsto \frac{1}{-\ln(z)}\right)^{\log} = \zeta \mapsto \ln(\zeta)$.\\
	
	As for transit maps between transversals on the same orbit, they are holomorphic maps with zero constants, by flow box they are even analytic functions with nonzero positive derivative let $f(z) = \alpha z + \sum_{q \geq 0}a_{q}z^{q + 2}$, $\alpha > 0$:
	
	\[f^{\log}(\zeta) = -\ln\left(\alpha e^{-\zeta} + \sum_{q \geq 0}a_{q}e^{-(q + 2)\zeta}\right) = \zeta - \ln(\alpha) - \ln\left(1 + \sum_{q \geq 0}\frac{a_{q}}{\alpha}e^{-(q + 1)\zeta}\right).\]
	
	Thus $f^{\log}$ has an expansion in terms of $e^{-k\zeta}$ with constant coefficients. Let $R$ be the radius of convergence of $f$, then this expansion is valid for $a$ large enough in:
	
	\[\mathbb{C}_{a}^{+} \coloneqq \{\zeta \in \mathbb{C} \mid \re(z) \geq a\}.\]
	
	Indeed, note that going into the logarithmic chart turns small discs around the origin into half planes of the form $\mathbb{C}_{a}^{+}$.\\
	
	As for the map $h: z \mapsto z^{\alpha}$:
	
	\[h^{\log}(\zeta) = \alpha \zeta.\]
	
	This means that going through an equilibrium from the hyperbolic side is given in the logarithmic chart first by multiplication by $k$, then by the exponential function, and for going through an equilibrium from the center side it is division by $k$ composed with the logarithm.\\
	
	Based on these two calculations we introduce two classes of functions:
	
	\begin{Definition}
		Let us define $\Aff$ to be the set of all maps of the form:
		
		\[\zeta \mapsto \alpha \zeta + \beta\]
		
		with $\alpha, \beta \in \mathbb{R}$ and $\alpha > 0$.
	\end{Definition}
	
	\begin{Definition}
		Let us define $\mathcal{H}$ to be the set of analytic functions $f$ for which there exists some series:
		
		\[\zeta + \sum_{q \geq 0}b_{q}e^{-c_{q}\zeta}\]
		
		$b_{k}, c_{k} \in \mathbb{R}$, $c_{k} > 0$ and strictly increasing to $+\infty$, such that for some $a \geq 0$, for any $\zeta \in \mathbb{C}^{+}_{a}$ and for any $N > 0$ and any $\epsilon \in (0, c_{N + 1} - c_{N})$ we have:
		
		\[\left|f(\zeta) - \left(\id + \sum_{q = 0}^{N}b_{k}e^{-c_{q}\zeta}\right)\right| \leq e^{-(c_{N} + \epsilon)\re(\zeta)}.\]
		
		We call the series above the Dulac series of $f$. (It is in fact true as we will elaborate on later, that this series is unique if it exists.)
	\end{Definition}
	
	Obviously we need these definitions because any Dulac map can be decomposed into these maps together with $\exp$ and $\ln$.
	
	\section{Strategy of proof}\label{SectFormal}
	
	Let us take a closer look at the structure of these return maps $\Delta$. We made the condition that the transit map of any equilibrium we pass is either of exponential type or of logarithmic type and that they alternate.\\
	
	We also made the condition that these types alternate, there are an even amount of equilibria and we assumed the first type we encountered was of exponential type. Let $\langle . \rangle$ be the group of germs of real analytic functions generated by the functions between the brackets under composition, then:
	
	\[\Delta \in \langle \Aff, \mathcal{H}, \ln\circ \Aff \circ \mathcal{H} \circ \exp\rangle.\]
	
	We can simplify this a bit further by noting that:
	
	\[\ln\circ \Aff \circ \mathcal{H} \circ \exp = \ln\circ \Aff \circ \exp \circ \ln \circ \mathcal{H} \circ \exp.\]
	
	And let $a \in \Aff, a(x) = \alpha x + \beta$, then:
	
	\[\ln(a(e^{\zeta})) = \ln(\alpha e^{\zeta} + \beta) = \zeta + \ln(\alpha) + \ln\left(1 + \frac{\beta}{\alpha} e^{-\zeta}\right)\]
	
	thus $\ln\circ \Aff \circ \exp \in \mathcal{H}$. Thus we have:
	
	\begin{Theorem}[Structural Theorem, specific case]
		Let $A$ stand for conjugation by $\exp$, i.e.:
		
		\[A(h) = \ln \circ h \circ \exp.\]
		
		All Dulac maps of the polycycles we study are contained in:
		
		\[\langle \Aff, \mathcal{H}, A\mathcal{H}\rangle.\]
	\end{Theorem}
	
	\begin{Remark}
		This is a specific case of \cite[\S 1.3 Proposition 3]{ilyashenkoFiniteness}.
	\end{Remark}
	
	The idea of \cite{ilyashenkoFiniteness} is this: $\Aff$ is far from identity, but $\mathcal{H}$ is close to identity and $A(\mathcal{H})$ is even closer to identity. More concretely, any element of $\mathcal{H}$ can be written as $\id + \phi$ with $\phi$ exponentially small and any element of $A(\mathcal{H})$ can be written as $\id + \psi \circ \exp$ with $\psi$ exponentially small, simply by applying $A$ to the Dulac series of an element in $\mathcal{H}$.\\
	
	The set of elements $\phi$ and $\psi$ we are willing to allow will be called functional cochains, the set of $\phi$ will be $FC^{0}$, functional cochains of level $0$ and the set of $\psi$ will be $FC^{1}$, functional cochains of level $1$.
	
	\begin{Remark}
		In \cite{ilyashenkoFiniteness} there is a further distinction between $FC^{n}_{0}$ and $FC^{n}_{1}$ which has to do with general semihyperbolic saddles/saddle nodes not being real analytically orbitally equivalent to the system we have here. This messes with domains as soon as you encounter $n + 1$ more exponential type saddles than logarithmic type at some point along the Dulac map. In some sense the $FC^{n}_{0}$ are the ones with large domain and $FC^{n}_{1}$ are the ones with small domain.
	\end{Remark}
	
	The idea is this, it is possible to expand all these compositions into sums using a lot of Taylor expansions and if we had a sum such as:
	
	\[a + \phi + \psi\circ \exp\]
	
	then $a$ could not interfere with $\phi$, $\phi$ could not interfere with $\psi$, because they are all of intrisically different sizes, so we could look in this sum for leading terms because we have separated the `levels of asymptotics' which we talked about in the introduction. So what we are looking for is a way to expand these compositions into sums `split by levels of asymptotics' indeed roughly corresponding to how $\Aff$ is far from the identity, $\mathcal{H}$ is close to the identity and $A\mathcal{H}$ is even much closer. This is the point of the so-called additive decomposition theorem \cite[p 73]{ilyashenkoFiniteness}.\\
	
	For this idea there are five crucial properties we want $FC^{0}$ to satisfy:
	
	\begin{enumerate}
		\item{$\mathcal{H} \subseteq \id + FC^{0}$.}
		\item{Let $a \in \Aff$, then:
		
		\[a \circ (\id + FC^{0}) \circ a^{-1} \subseteq \id + FC^{0}.\]}
		\item{Let $a \in \Aff$, then:
		
		\[a \circ (\id + FC^{0}) \subseteq a + FC^{0}.\]}
		\item{$(\id + FC^{0})$ forms a group under composition.}
		\item{Let $\phi \in FC^{0}$ be nonzero then there exists some $\lambda > 0$ such that for $x \in \mathbb{R}$ large enough:
		
		\[|\phi(x)| \geq e^{-\lambda x}.\]}
	\end{enumerate}
	
	The $FC^{0}$ proposed in \cite{ilyashenkoFiniteness} simplifies down to the following:
	
	\begin{Definition}
		The set $FC^{0}$ is the set of analytic functions satisfying the same conditions as $\mathcal{H}$ with Dulac series of the form:
		
		\[\sum_{q \geq 0}b_{q}e^{-c_{q}\zeta}.\]
	\end{Definition}
	
	\begin{Remark}
		In \cite{ilyashenkoFiniteness} there could also be hyperbolic equilibria in the polycycle, this complicates both the domain and the asymptotics, we will not get into the domain, as for the asymptotics, the $b_{q}$ are replaced by (complex) polynomials in $\zeta$.
	\end{Remark}
	
	\begin{Remark}\label{RemFC+}
		We will for ease of notation restrict ourselves to the $FC^{n}$ which are exponentially decreasing themselves. This is not the case in \cite{ilyashenkoFiniteness}, but wherever that is relevant we will instead choose to write out very concretely what happens in \cite{ilyashenkoFiniteness}.\\
		
		More specifically if we were to follow the notation of \cite{ilyashenkoFiniteness} then $FC^{0}$ would be the above set without the restriction that $c_{0} > 0$ and the set we describe here would be $FC^{0}_{+}$.
	\end{Remark}
	
	By some routine calculation it is possible to verify all five conditions, upon being given the following (see e.g. \cite[\S 3.1 Corollary 1]{ilyashenkoFiniteness})
	
	\begin{Theorem}\label{ThmPL}
		Let $f$ be a bounded holomorphic function on $\mathbb{C}^{+}_{0}$ decreasing faster than any exponential on the real axis, then $f$ is identically zero.
	\end{Theorem}
	
	The argument essentially being the following, suppose $\phi \in FC^{0}$ has a nonzero series, then by leading term we get the exponential lower bound we want. Suppose this asymptotic series is zero, then $\phi$ is smaller than any exponential on the real axis and bounded, thus identically zero.\\
	
	The way we get the promised sum out of these properties is the following, suppose that we have a particularly nice $\Delta$, in the sense that it is inside $\langle \Aff, \mathcal{H}\rangle$, then using these properties we can conjugate elements of $\Aff$ past elements of $(\id + FC^{0}) \supset \mathcal{H}$ and expand in order to write:
	
	\[\Delta \in \Aff + FC^{0}.\]
	
	\begin{Remark}
		This is a version of the central Theorem that \cite{ilyashenkoFiniteness} works towards, it is a consequence of $ADT_{1}$ as stated in \cite[p 73]{ilyashenkoFiniteness}
	\end{Remark}
	
	Then $\Delta - \id$ is in the same set and either $\Aff - \id$ is nonzero, giving a zero free region because an element of $FC^{0}$ is much smaller, or $FC^{0}$ is nonzero and $\Aff - \id$ is zero, giving by the lower bound of $FC^{0}$ a zero free region, or $\Delta - \id \equiv 0$, thus it having no isolated zeroes.\\
	
	This is essentially a leading term argument dating back in a flawed form to the arguments of Dulac.\\
	
	We now want to continue this story of rearranging by conjugation and writing out as a sum, where we replace $A\mathcal{H}$ by $\id + FC^{1} \circ \exp$, this is exactly what we do with elements of $(\id + FC^{0})$, but with $\Aff$ the story is slightly different.\\
	
	Taking an element $f \in A\mathcal{H}$ we can quickly calculate the series it has:
	
	\begin{align*}
		f(\zeta) & = \ln\left(e^{\zeta} + \sum_{q \geq 0} b_{q}e^{-qe^{\zeta}}\right)\\
		& = \zeta + \ln\left(1 + \sum_{q \geq 0} b_{q}e^{-\zeta}e^{-qe^{\zeta}}\right)\\
		&= \zeta + \sum_{q \geq 0}P_{q}(e^{-\zeta})e^{-c_{q}e^{\zeta}}\\
	\end{align*}
	
	with $P_{k}$ being polynomials, $c_{k} > 0$, strictly increasing to $+\infty$. The last step going by Taylor expansion, which is fine as soon as we restrict the domain to some $\ln(\mathbb{C}_{a}^{+})$. Then $f - \id$ is roughly of order $e^{-\lambda \exp(\zeta)}$ for some $\lambda > 0$. Let $a \in \Aff$, $a(\zeta) = \alpha \zeta + \beta$, then:
	
	\[(a^{-1}\circ f \circ a)(\zeta) = \zeta + \sum_{k \geq 0}\tilde{P}_{k}(e^{-\alpha\zeta})e^{-\tilde{c}_{k}e^{\alpha\zeta}}.\]
	
	So if $\alpha \neq 1$, then this is roughly of size $e^{-\lambda e^{\alpha \zeta}}$ for some $\lambda > 0$, which is of essentially different size than $e^{-\mu e^{\zeta}}$ for all $\mu > 0$. So should we want as in \cite{ilyashenkoFiniteness} that $FC^{1}$ has exponential upper and lower bound, then:
	
	\[a^{-1}\circ(\id + FC^{1}\circ \exp) \circ a \not\subseteq \id + FC^{1}\circ \exp.\]
	
	At best we can expect the following, let $m_{\alpha}$ be the map $\zeta \mapsto \alpha \zeta$, then we can ask for:
	
	\[a^{-1}\circ(\id + FC^{1}\circ \exp) \circ a \subseteq \id + FC^{1}\circ \exp \circ m_{\alpha}.\]
	
	\begin{Remark}
		This choice to hold on to the exponential lower and upper bound instead of making one big $FC^{1}$ is something intimately related to the quasianalyticity properties.\\
		
		The reason the quasianalyticity property in \cite{ilyashenkoFiniteness} works is that after extending using small Stokes phenomena we get control of the functions in $FC^{n}$ up to a a domain of `essentially the same size as $\mathbb{C}^{+}_{a}$', in the sense that Theorem \ref{ThmPL} holds for the domain. The point is then that because the Stokes phenomena are so small, essentially `one level of asymptotics lower' then any element of $FC^{n}$ smaller than any exponential on the real axis will be identically zero on the real axis.\\
		
		In the case we were to let go of the exponential lower bound and merge all the $FC^{1}$ we would have functions of size $e^{-x^{\alpha}}$ for all $\alpha > 0$, in particular for $\alpha > 1$ we would have things smaller than any exponential, so we would ruin all quasianalyticity arguments.
	\end{Remark}
	
	Then the minimal list of conditions on $FC^{1}$ is the following, essentially taken from the conditions needed in \cite[\S 1.11]{ilyashenkoFiniteness}:
	
	\begin{enumerate}
		\item{$A\mathcal{H} \subseteq \id + FC^{1} \circ \exp$.}
		\item{Let $\alpha > 0$ and $m_{\alpha}$ be multiplication with $\alpha$. Then:
		
		\begin{enumerate}
			\item{$\id + FC^{1}\circ \exp \circ m_{\alpha}$ is a group under composition.}
			\item{Let $f \in \id + FC^{0}$, then:
			
			\[f^{-1} \circ (\id + FC^{1}\circ \exp \circ m_{\alpha}) \circ f \subseteq \id + FC^{1}\circ \exp \circ m_{\alpha}.\]}
			\item{Let $a \in \Aff$:
			
			\[a \circ (\id + FC^{1}\circ \exp \circ m_{\alpha}) \subseteq a + FC^{1} \circ \exp \circ m_{\alpha}.\]}
			\item{Let $\phi \in FC^{0}$:
			
			\[\phi \circ (\id + FC^{1}\circ \exp \circ m_{\alpha}) \subseteq \phi + FC^{1} \circ \exp \circ m_{\alpha}.\]}
		\end{enumerate}}
		\item{Let $\alpha_{1} > \alpha_{2} > 0$ and let $m_{\alpha_{i}}$ be multiplication with $\alpha_{i}$. Then:
		
		\begin{enumerate}
			\item{Let $f \in \id + FC^{1} \circ \exp \circ m_{\alpha_{2}}$, then:
			
			\[f^{-1} \circ (\id + FC^{1} \circ \exp \circ m_{\alpha_{1}}) \circ f \subseteq \id + FC^{1} \circ \exp \circ m_{\alpha_{1}}.\]}
			\item{Let $g \in FC^{1} \circ \exp \circ m_{\alpha_{2}}$, then:
			
			\[g \circ (\id + FC^{1} \circ \exp \circ m_{\alpha_{1}}) \subseteq g + FC^{1} \circ \exp \circ m_{\alpha_{1}}.\]}
			\item{Let $i, j = 1, 2$, let $a \in \Aff$, $a(\zeta) = \alpha_{i} \zeta + \beta$, then:
			
			\[a^{-1} \circ (\id + FC^{1} \circ \exp \circ m_{\alpha_{j}}) \circ a \subseteq \id + FC^{1} \circ \exp \circ m_{\alpha_{j}}\circ m_{\alpha_{i}}.\]}
		\end{enumerate}}
		\item{Let $\psi \in FC^{1}$ be nonzero then there exists some $\lambda > 0$ such that for $x \in \mathbb{R}$ large enough:
			
			\[|\psi(x)| \geq e^{-\lambda x}.\]}
	\end{enumerate}
	
	The idea is the same, you conjugate past and expand until you can write every Dulac map $\Delta$ as:
	
	\[\Delta \in \Aff + FC^{0} + \sum_{p = 1}^{q}FC^{1} \circ \exp \circ m_{\alpha_{p}}\]
	
	each of the $\alpha_{p}$ different and dependent on $\Delta$ (obviously $q$ also depends on $\Delta$).\\ 
	
	\begin{Remark}
		This is the full importance of $ADT_{1}$ in \cite[p 73]{ilyashenkoFiniteness} in this context.
	\end{Remark}
	
	Note that by our earlier remark that for $\alpha \neq 1$ a function of size $e^{-\lambda\exp^{\alpha \zeta}}$,  is of essentially different size than $e^{-\mu\exp^{\mu\zeta}}$, the elements in different $FC^{1} \circ \exp \circ m_{\alpha_{p}}$ can not have any overlap. So this gives us a leading term for $\Delta - \id$ or $\Delta \equiv \id$.\\
	
	The problem in \cite{ilyashenkoFiniteness} lies in its proof of property $3$, more specifically $2b$ significantly complicates the asymptotics of $FC^{1}$ and the proof that these asymptotics are ordered is problematic.
	
	\section{Sketch of definition of $FC^{1}$ and ordering of asympotics}\label{SectProblem}
	
	In \cite[\S 1.7]{ilyashenkoFiniteness}, $FC^{1}$ is defined as the union of $FC^{1, p}$ which are defined by inductive process as follows (again some simplification are made appropriately for the context):\\
	
	First define the set of exponents $E^{1}$ to be the set of all functions $\mathbf{e}$ which on some $\mathbb{C}^{+}_{a}$ admits asymptotics of the form:
	
	\[\mathbf{e}(\zeta) = \nu(\mathbf{e})\cdot e^{\zeta} + \sum_{q \geq 0}b_{q}e^{(1 - c_{q})\zeta}\]
	
	again, all $c_{q} > 0$, strictly increasing, going to $+\infty$. We call $\nu(\mathbf{e}) \in \mathbb{R}$ the principal exponent of $\mathbf{e}$.\\
	
	We define $K_{1, 0}$ to be sums of $FC^{0}$ multiplied by some exponential $e^{\lambda \zeta}$, $\lambda > 0$ (i.e. the $FC^{0}$ in the notation of \cite{ilyashenkoFiniteness} as explained in Remark \ref{RemFC+}).\\
	
	We proceed by induction, suppose $K_{1, p}$ defined, then we define $FC^{1, p}$ to be the set of all analytic functions $\psi$ with $\psi \circ \exp$ admitting double exponentially accurate asymptotics on some $\ln(\mathbb{C}^{+}_{a})$ of the form:
	
	\[\sum_{q \geq 0} k_{q}e^{\mathbf{e}_{q}}\]
	
	with $k_{q} \in K_{1, p}$ and $\mathbf{e}_{q} \in E^{1}$. I.e. for a given $\lambda > 0$ any initial segment after a given point approximates the element in $FC^{1, p}$ closer than $e^{-\lambda \re(e^{\zeta})}$ on the real axis. We also demand that:
	
	\[\lim_{q \to \infty}\nu(\mathbf{e}_{q}) = -\infty.\]
	
	It is worth noting that multiple $\mathbf{e}_{q}$ can have the same principle exponent, we just assume them ordered in some (not necessarily strictly) descending order.
	
	Then $K_{1, p + 1}$ is a sum of the following form:
	
	\[\phi + \sum_{q = 1}^{N} \psi_{q}\circ \exp \circ m_{\alpha_{q}}\]
	
	with $\phi \in K_{1, 0}$, $\psi_{q} \in FC^{1, p}$, each of the $m_{\alpha_{q}}$ multiplication with $\alpha_{q}$ with $0 < \alpha_{q} < 1$. Then as said before $FC^{1}$ is the union of all the $FC^{1, p}$.
	
	\begin{Remark}
		An expression of the form:
		
		\[\sum_{q \geq 0} k_{q}e^{\mathbf{e}_{q}}\]
		
		is called a $\STAR$-series. This notion is less well behaved than one might expect, in particular it is possible to have a nontrivial finite sum of the form above which is nontheless zero, indeed for example $e^{-e^{\zeta}} - 1$ can be expanded into Taylor series to give an element of $FC^{0}$.\\
		
		The focus is nonetheless on finite sums which are nonzero.
	\end{Remark}
	
	Let us unwrap the proof that a nonzero element $\psi$ of $FC^{1}$ has an exponential lower bound as done in \cite{ilyashenkoFiniteness}.\\
	
	The first step is to use asymptotics to create a dichotomy: either $\psi$ is smaller than any exponential on the real axis, or it has the lower bound (this is the problematic part) and then we prove that if $\psi$ is smaller than any exponential on the real axis it is identically zero, see \cite[pp 74--75]{ilyashenkoFiniteness}, \cite[Chapter III]{ilyashenkoFiniteness} is dedicated to the statement that if $\psi$ is smaller than any exponential, it is identically zero, but in this case it is the same Theorem as before.\\
	
	The way to create this dichotomy is essentially found in \cite[\S 4.10 G]{ilyashenkoFiniteness}. The point is to prove by induction on $p$ that any nonzero finite sum:
	
	\[\sum_{q = 1}^{N} k_{q}e^{\mathbf{e}_{q}}\]
	
	has a (sharp) exponential lower bound of exactly the same type demanded with $k_{q} \in K_{1, p}$. So if we have a series approximating an element of $FC^{1, p} \circ \exp$, either eventually this lower bound is higher than the accuracy by which it is supposed to approximate this element of $FC^{1, p} \circ \exp$, implying this lower bound for the element of $FC^{1, p}$ or the element of $FC^{1, p}$ is smaller than any exponential on the real axis, making it zero.\\
	
	Here we proceed by induction on $N$, i.e. the amount of terms in such a sum. The core of this is a divide and differentiate argument. The case $N = 1$ is easy because an element of $K_{1, p}$ is intrinsically smaller than an element of $\exp(E^{1})$ (remembering that elements of $FC^{1}$ are supposed to be exponentially small).\\
	
	Then by ordering already proven one can assume that $k_{1}e^{\mathbf{e}_{1}}$ is the largest and one can consider:
	
	\[S(x) \coloneqq 1 + \sum_{q = 2}^{N}\frac{k_{q}}{k_{1}}e^{\mathbf{e}_{q} - \mathbf{e}_{1}}.\]
	
	Then assuming this has a limit as $x$ goes to infinity (this is the first problematic part) we can consider this limit which has to be $\leq 1$ in absolute value, if it is nonzero we are done because this is roughly a nonzero constant times something already having the lower bound we want.\\
	
	Suppose it does go to zero, then:
	
	\[S(x) = \int_{x}^{\infty}S'(y)dy.\]
	
	And $S'(y)$ is the same kind of finite sum with $N - 1$ terms (this is also problematic), so a lower bound for $S'(y)$ implies by some elementary integral estimates a lower bound for $S$ and thus one for the original sum we started with.\\
	
	\section{The problem and a counterexample}\label{SectCounterExample}
	
	The problem here is that both of the passages where we say it is problematic make the following calculation:
	
	\[\frac{d}{d\zeta}\frac{k_{q}}{k_{1}}e^{\mathbf{e}_{q} - \mathbf{e}_{1}} = \frac{k_{q}'k_{1} - k_{1}'k_{q} + k_{q}k_{1}(\mathbf{e}_{q}' - \mathbf{e}_{1}')}{k_{1}^{2}}e^{\mathbf{e}_{q} - \mathbf{e}_{1}}.\]
	
	And then claim that $k_{q}'k_{1} - k_{1}'k_{q} + k_{q}k_{1}(\mathbf{e}_{q}' - \mathbf{e}_{1}')$ is in $K_{1, p}$, see \cite[p198 last line]{ilyashenkoFiniteness}, calculation slightly corrected (technically \cite{ilyashenkoFiniteness} also claims falsely in the same line that $k_{1}^{2}$ is in $K_{1, p}$ but that part is not essential to the argument).\\
	
	\begin{Remark}
		To make this explicit, in \cite[p198 last line]{ilyashenkoFiniteness} it talks about $\mathscr{K}^{m, r}_{\mathbb{R}}$. We have:
		
		\[K_{1, p} \subset \mathscr{K}^{1, p}_{\mathbb{R}}\]
		
		moreover the same arguments we will employ will indeed also work to show that in the counterexample we have:
		
		 \[k_{q}'k_{1} - k_{1}'k_{q} + k_{q}k_{1}(\mathbf{e}_{q}' - \mathbf{e}_{1}') \notin \mathscr{K}^{1, p}_{\mathbb{R}}\]
	\end{Remark}
	
	This is true for $p = 0$. We will show that this can not be the case for $p = 1$ using maps coming from concrete polycycles in order to show that this difficulty is incircumventable in this approach. What we mean with this is the following: Our example will show that:
	
	\begin{enumerate}
		\item{This $FC^{1, p}$ construction is inevitable, i.e. we will show that we naturally get an element of $FC^{1, 1}$.}
		\item{For elements $k_{2}, k_{1}$ of $K_{1, 1}$ coming from this example we will show that:
		
		\[k_{2}'k_{1} - k_{1}'k_{2} \notin K_{1, 1},\]
		
		giving a counterexample to the claims made in proving the ordering of $FC^{1}$ and thus the existence of leading terms, or the lower bound needed.}
	\end{enumerate}
	
	The crux of the entire thing is that you forcibly have that $FC^{1, p}\circ \exp \cdot FC^{1, 0}\circ \exp \circ m_{\alpha} \subseteq FC^{1, p + 1}$ with $\alpha < 1$ because any approximation you try to give the element of $FC^{1, 0}\circ \exp \circ m_{\alpha}$ will be off by something larger than $e^{-x}$ so you have no choice but to use the entire element of $FC^{1, 0}\circ \exp \circ m_{\alpha}$ in your coefficients, raising the $p$ by one.\\
	
	Let us get to the counterexample. We take a polycycle with $4$ equilibria $E_{1}, ..., E_{4}$ in order in forward time. Let us take the section $\Sigma$ around which to take a return map. Let $\Sigma_{i, 1}$ and $\Sigma_{i, 2}$ be two sections around $E_{i}$ as in section 2 with a transit map going from $\Sigma_{i, 1}$ to $\Sigma_{i, 2}$. 
	
	\begin{figure}[H]
		\centering
		\includegraphics[scale = 1]{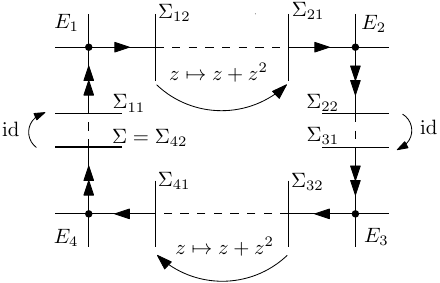}
		\caption{Counterexample notation}
	\end{figure}
	
	Let us prescribe all the transit maps by gluing. What this means is that for each equilibrium we have a concrete system on $\mathbb{R}^{2}$ with the equilibrium, say at the origin, with the correct transit map. Then it is possible to use flow box theorem to glue both the manifolds and the vector field at the same time and we can use a biholomorphism in the direction perpendicular to the orbits of the flow box to prescribe the transit maps between different equilibria:
	
	\begin{enumerate}
		\item{$\Sigma \to \Sigma_{1, 1}: \id.$}
		\item{$\Sigma_{1, 1} \to \Sigma_{1, 2}: z \mapsto e^{-\frac{1}{z^{2}}}.$}
		\item{$\Sigma_{1, 2} \to \Sigma_{2, 1}: z \mapsto z + z^{2}.$}
		\item{$\Sigma_{2, 1} \to \Sigma_{2, 2}: z \mapsto \left(\frac{1}{-\ln(z)}\right)^{\frac{1}{2}}.$}
		\item{$\Sigma_{2, 2} \to \Sigma_{3, 1}: \id$.}
		\item{$\Sigma_{3, 1} \to \Sigma_{3, 2}: z \mapsto e^{-\frac{1}{z}}$.}
		\item{$\Sigma_{3, 2} \to \Sigma_{4, 1}: z \mapsto z + z^{2}$.}
		\item{$\Sigma_{4, 1} \to \Sigma_{4, 2}: z \mapsto \frac{1}{-\ln(z)}$.}
		\item{$\Sigma_{4, 2} \to \Sigma: \id$.}
	\end{enumerate}
	
	Let us try to convert these maps into the logarithmic chart. Let us first consider the map $f(z) = z + z^{2}$:
	
	\[f^{\log}(\zeta) = \zeta - \ln(1 + e^{-\zeta}) = \zeta - \sum_{q > 0}(-1)^{q}\frac{e^{-q\zeta}}{q}\]
	
	\[A(f^{\log})(\zeta) = \zeta - \ln\left(1 - \sum_{q > 0}(-1)^{q}e^{-\zeta}\frac{e^{-qe^{\zeta}}}{q}\right) = \zeta + \sum_{r > 0}\frac{\left(\sum_{q > 0}(-1)^{q}e^{-\zeta}\frac{e^{-qe^{\zeta}}}{q}\right)^{r}}{r}\]
	
	\[(m_{2}^{-1} \circ A(f^{\log}) \circ m_{2})(\zeta) = \zeta + \frac{1}{2}\sum_{r > 0}\frac{\left(\sum_{q > 0}(-1)^{q}e^{-2\zeta}\frac{e^{-qe^{2\zeta}}}{q}\right)^{r}}{r}.\]
	
	The return map is then:
	
	\[A(f^{\log}) \circ (m_{2}^{-1} \circ A(f^{\log}) \circ m_{2}).\]
	
	Note that this could only be identity if $A(f^{\log})\circ m_{2} = m_{2} \circ A(f^{\log})^{-1}$ which is clearly not the case, they are of size $2\zeta + e^{-\lambda e^{2\zeta}}$ and $2\zeta + e^{-\mu e^{\zeta}}$ respectively for some $\lambda, \mu > 0$. We now promised that this example would illustrate two things, so let us prove this.\\
	
	\underline{\textbf{Part 1: the $FC^{1, p}$ construction is necessary}}\\
	
	The explicit form of the return map means that by property $2b$ of $FC^{1}$:
	
	\begin{multline*}
		\left(\sum_{r > 0}\frac{\left(\sum_{q > 0}(-1)^{q}e^{-\zeta}\frac{e^{-qe^{\zeta}}}{q}\right)^{r}}{r}\right)\circ \left(\zeta + \frac{1}{2}\sum_{r > 0}\frac{\left(\sum_{q > 0}(-1)^{q}e^{-2\zeta}\frac{e^{-qe^{2\zeta}}}{q}\right)^{r}}{r}\right)-\\ - \sum_{r > 0}\frac{\left(\sum_{q > 0}(-1)^{q}e^{-\zeta}\frac{e^{-qe^{\zeta}}}{q}\right)^{r}}{r} \in FC^{1} \circ \exp \circ m_{2}.
	\end{multline*}
	
	It is easy to calculate this out using Taylor expansion:
	
	\[\sum_{s \geq 1}\frac{1}{s!}\frac{d^{s}\sum_{r > 0}\frac{\left(\sum_{q > 0}(-1)^{q}e^{-\zeta}\frac{e^{-qe^{\zeta}}}{q}\right)^{r}}{r}}{d\zeta^{s}}\left(\frac{1}{2}\sum_{r > 0}\frac{\left(\sum_{q > 0}(-1)^{q}e^{-2\zeta}\frac{e^{-qe^{2\zeta}}}{q}\right)^{r}}{r}\right)^{s}\]
	
	has to be in $FC^{1} \circ \exp \circ m_{2}$. Now the important part here is that the first factor, the entire thing being derived, has asymptotics of accuracy $e^{-\lambda e^{\zeta}}$ while the definition of asymptotics for $FC^{1, p}$ needs accuracy up to $e^{-\lambda e^{2\zeta}}$ so we have no choice but to take the entire thing as a coefficient. Put differently, what we have here is an element of $FC^{1, 1}\circ \exp \circ m_{2}$. So $FC^{1, 0}$ definitely does not cover everything. So in the approach of \cite{ilyashenkoFiniteness} it is actually necessary that for $k_{1}, k_{2} \in K_{1, 1}$:
	
	\begin{equation}\label{EqCounter}
		k_{2}'k_{1} - k_{1}'k_{2} \in K_{1, 1}
	\end{equation}
	
	(here we have taken the previous exponents equal to each other).\\
	
	\underline{\textbf{Part 2: the argument for ordering does not work: theoretical example}}\\
	
	So we want to prove that Equation \ref{EqCounter} does not hold in a case relevant to proving finiteness of limit cycles. We will get to a case actually related to polycycles later.\\ 
	
	Let us first illustrate the problem in a much simpler context removed from the problem of Dulac. Let:
	
	\[k_{1}(\zeta) = \frac{e^{-e^{\frac{1}{3}\zeta}}}{1 - e^{-e^{\frac{1}{3}\zeta}}} \in FC^{1, 0}\circ \exp \circ m_{\frac{1}{3}} \subseteq K_{1, 1},\]
	
	\[k_{2}(\zeta) = e^{-e^{\frac{1}{2}\zeta}} \in FC^{1, 0}\circ \exp \circ m_{\frac{1}{2}} \subseteq K_{1, 1}.\]
	
	In fact for these it also holds that $k_{2}'k_{1} - k_{1}'k_{2} \notin K_{1, 1}$ by the same issue we will talk about now but let us prove that $k_{1}k_{2} \notin K_{1, 1}$ to keep calculation to a minimum. It is worth noting that Equation \ref{EqCounter} would be very extraordinary if $K_{1, 1}$ were not at least a differential algebra, so in some sense it is also expected that $k_{1}k_{2} \in K_{1, 1}$:
	
	\[k_{1}(\zeta)k_{2}(\zeta) = \frac{e^{-e^{\frac{1}{3}\zeta}}}{1 - e^{-e^{\frac{1}{3}\zeta}}}e^{-e^{\frac{1}{2}\zeta}}.\]
	
	Suppose by contradiction it were in $K_{1, 1}$, thus it could be written as:
	
	\[\phi + \sum_{q = 1}^{N} \psi_{q}\circ \exp \circ m_{\alpha_{q}}\]
	
	with $\phi \in K_{1, 0}$, $\psi_{q} \in FC^{1, p}$, each of the $m_{\alpha_{q}}$ multiplication with $\alpha_{q}$ with $0 < \alpha_{q} < 1$. Then clearly because $k_{1}k_{2}$ has both an upper and lower bound of the form $e^{-\lambda e^{\frac{1}{2}\zeta}}$ we would have to have that this is a sum of $FC^{1, 0}\circ \exp \circ m_{\alpha_{q}}$ with $\alpha_{q} \leq \frac{1}{2}$.\\
	
	So by construction of $K_{1, 0}$ there has to be some element of $FC^{1, 0} \circ \exp \circ m_{\frac{1}{2}}$ which approximates $k_{1}k_{2}$ with accuracy $e^{-\lambda e^{\frac{1}{2}\zeta}}$ for all $\lambda$. Now an element of $FC^{1, 0} \circ \exp$ looks like:
	
	\[\sum_{q = 1}^{N} h_{q}e^{\mathbf{e}_{q}}.\]
	
	With $h_{q}$ having some exponential lower and upper bound, in particular we have:
	
	\[k_{1}(\zeta)k_{2}(\zeta) = \frac{e^{-e^{\frac{1}{3}\zeta}}}{1 - e^{-e^{\frac{1}{3}\zeta}}}e^{-e^{\frac{1}{2}\zeta}} = \left(\sum_{q\geq 1}e^{-qe^{\frac{1}{3}\zeta}}\right)e^{-e^{\frac{1}{2}\zeta}} = \sum_{q\geq 1}e^{-qe^{\frac{1}{3}\zeta}-e^{\frac{1}{2}\zeta}}.\]
	
	So a finite sum of $h_{q}e^{\mathbf{e}_{q}}$ has to approximate:
	
	\[\sum_{q\geq 1}e^{-qe^{\frac{2}{3}\zeta}-e^{\zeta}}\]
	
	up to accuracy $e^{-2e^{\zeta}}$.\\ 
	
	\underline{\textbf{Diversion: Some normal forms}}\\
	
	Let us introduce some normal forms for $he^{\mathbf{e}}$, $h \in K_{1, 0}$, $\mathbf{e} \in E^{1}$ in order to definitively show that this is impossible. Remember that:
	
	\[\mathbf{e}(\zeta) = \nu(\mathbf{e})\cdot e^{\zeta} + \sum_{q \geq 0}b_{q}e^{(1 - c_{q})\zeta}.\]
	
	So it is possible to split up a $\mathbf{e} \in E^{1}$ as $\mathbf{e}_{\Larget} + \mathbf{e}_{\Small}$ where $\mathbf{e}_{\Small}$ contains all terms with $c_{q} \geq 1$ and the important remark here is that by Taylor series:
	
	\[e^{\mathbf{e}_{\Small}} \in K_{1, 0}.\]
	
	So let $E^{1}_{\Larget}$ be the subset of all $E^{1}$ where the $c_{q}$ are all $< 1$ then any finite sum of $he^{\mathbf{e}}$, $h \in K_{1, 0}$, $\mathbf{e} \in E^{1}$ can be rewritten into a finite sum with $\mathbf{e} \in E^{1}_{\Larget}$. The important thing is the following, because of the exponential lower and upper bound of $K_{1, 0}$, if we have $h_{1}, h_{2} \in K_{1, 0}$ and $\mathbf{e}_{1} < \mathbf{e}_{2} \in E^{1}_{\Larget}$ then as we go to $+\infty$ on the real axis:
	
	\[\frac{h_{1}e^{\mathbf{e}_{1}}}{h_{2}e^{\mathbf{e}_{2}}} \to 0,\]
	
	because it behaves like some $e^{\lambda x}e^{-\mu e^{\alpha\zeta}}$ with $0 < \alpha < 1$ and of course the double exponential will win out.\\
	
	So we have rewritten this to a form with proper leading terms.\\
	
	\underline{\textbf{Back to part 2: the theoretical example}}\\
	
	From the normal forms it is clear that the only STAR-series one could choose to approximate $k_{1}k_{2}$ with coefficients in $K_{1, 0}$ is just:
	
	\[\sum_{q\geq 1}e^{-qe^{\frac{2}{3}\zeta}-e^{\zeta}}.\]
	
	But here the generalized exponent of $-qe^{\frac{2}{3}\zeta}-e^{\zeta} \in E^{1}_{\Larget}$ does not go to $-\infty$, in fact it stays $1$, thus it is not valid as a series expansion, thus $k_{1}k_{2} \notin K_{1, 0}$.\\
	
	\underline{\textbf{Part 3: the argument for ordering does not work: practical example}}\\
	
	Let us now go back to the example in Part 1 and prove that we can find some $k_{1}, k_{2}$ relevant to the example such that Equation \ref{EqCounter} does not hold.\\
	
	Let $k_{1} \in FC^{1, 0} \circ \exp \circ m_{\frac{1}{2}} \subseteq K_{1, 1}$ be:
	
	\[\sum_{r > 0}\frac{\left(\sum_{q > 0}(-1)^{q}e^{-\frac{1}{2}\zeta}\frac{e^{-qe^{\frac{1}{2}\zeta}}}{q}\right)^{r}}{r}\]
	
	so $\left(A(f^{\log}) - \id\right)\circ m_{\frac{1}{2}}$.\\
	
	Let $k_{2} \in FC^{1, 0} \circ \exp \circ m_{\frac{1}{3}} \subseteq K_{1, 1}$ be:
	
	\[\sum_{r > 0}\frac{\left(\sum_{q > 0}(-1)^{q}e^{-\frac{1}{3}\zeta}\frac{e^{-qe^{\frac{1}{3}\zeta}}}{q}\right)^{r}}{r}\]
	
	so $\left(A(f^{\log}) - \id\right)\circ m_{\frac{1}{3}}$. Then $k_{1}'$ is equal to:
	
	\[\left(\sum_{q > 0}-\frac{1}{2}(-1)^{q}e^{-\frac{1}{2}\zeta}\frac{e^{-qe^{\frac{1}{2}\zeta}}}{q}  -\frac{1}{2}(-1)^{q}e^{-qe^{\frac{1}{2}\zeta}}\right)\sum_{r > 0}\left(\sum_{q > 0}(-1)^{q}e^{-\frac{1}{2}\zeta}\frac{e^{-qe^{\frac{1}{2}\zeta}}}{q}\right)^{r - 1}\]
	
	or
	
	\[-\frac{1}{2}\left(\sum_{q > 0}(-1)^{q}e^{-qe^{\frac{1}{2}\zeta}}\left(\frac{e^{-\frac{1}{2}\zeta}}{q}  + 1\right)\right)\sum_{r > 0}\left(\sum_{q > 0}(-1)^{q}e^{-\frac{1}{2}\zeta}\frac{e^{-qe^{\frac{1}{2}\zeta}}}{q}\right)^{r - 1}.\]
	
	So $k_{2}'$ is equal to:
	
	\[-\frac{1}{3}\left(\sum_{q > 0}(-1)^{q}e^{-qe^{\frac{1}{3}\zeta}}\left(\frac{e^{-\frac{1}{3}\zeta}}{q}  + 1\right)\right)\sum_{r > 0}\left(\sum_{q > 0}(-1)^{q}e^{-\frac{1}{3}\zeta}\frac{e^{-qe^{\frac{1}{3}\zeta}}}{q}\right)^{r - 1}.\]
	
	It is then clear that if $k_{2}'k_{1} - k_{1}'k_{2} \in K_{1, 1}$ by the rough size being $e^{-\lambda e^{\frac{1}{2}\zeta}}$ the only option is that $k_{2}'k_{1} - k_{1}'k_{2} \in FC^{1, 0}\circ \exp \circ m_{\frac{1}{2}}$, but clearly you again by these normal forms can not write this as such.
	
	\section{Statement \ref{StatementPos}}\label{SectPosres}
	
	As promised we will talk about the case where every equilibrium is real analytically orbitally equivalent to:
	
	\[\begin{cases*}
		\dot{x} = x^{2}\\
		\dot{y} = -y.
	\end{cases*}\]
	
	In this case each of the transit maps is equal to $e^{-\frac{1}{z}}$. The crucial consequence is the following: we no longer just have that the Dulac maps in the logarithmic chart are contained in $\langle \Aff, \mathcal{H}, A\mathcal{H}\rangle$, we can also restrict to the case where every element in $\Aff$ has linear part $1$, i.e. is of the form $\zeta \mapsto \zeta + \beta$, after all the only multiplications, i.e. the only $\alpha$ come from the map $x \mapsto x^{k}$ which is no longer present. This avoids all this frustration with $m_{\alpha}, \alpha \neq 1$, allowing us to stay within $FC^{1, 0}$ and only use $K_{1, 0}$.\\ 
	
	The normal forms above work here to show ordering and even give proper leading terms. In a future article we hope to use these types of normal forms in a more general context.
	
	\section{Bibliography}
	
	\bibliography{mybib}
	\bibliographystyle{plain}
	
\end{document}